\documentclass[10pt]{amsart}
\usepackage{amsmath, amsthm, amssymb}
\usepackage{mathbbol}
\usepackage{txfonts}

\newtheorem{thm}{Theorem}[section]
\newcommand{\bt}{\begin{thm}}
\newcommand{\et}{\end{thm}}

\newtheorem{cor}[thm]{Corollary}
\newcommand{\bc}{\begin{cor}}
\newcommand{\ec}{\end{cor}}

\newtheorem{lem}[thm]{Lemma}
\newcommand{\bl}{\begin{lem}}
\newcommand{\el}{\end{lem}}

\newtheorem{prop}[thm]{Proposition}
\newcommand{\bp}{\begin{prop}}
\newcommand{\ep}{\end{prop}}

\newtheorem{defn}[thm]{Definition}
\newcommand{\bd}{\begin{defn}}      
\newcommand{\ed}{\end{defn}}

\newtheorem{rmrk}[thm]{Remark}
\newcommand{\br}{\begin{rmrk}}
\newcommand{\er}{\end{rmrk}}

\newtheorem{example}[thm]{Example}

\newcommand{\thmref}[1]{Theorem~\ref{#1}}
\newcommand{\secref}[1]{Section~\ref{#1}}
\newcommand{\lemref}[1]{Lemma~\ref{#1}}
\newcommand{\defref}[1]{Definition~\ref{#1}}

\newcommand{\N}{\mathbb{N}}

\newcommand{\R}{\mathbb{R}}
\newcommand{\Z}{\mathbb{Z}}

\newcommand{\hm}{{\mathcal H}}

\newcommand{\rstr}{\:\mbox{\rule{0.1ex}{1.2ex}\rule{1.1ex}{0.1ex}}\:}
\newcommand{\bdry}{\partial}

\newcommand{\spt}{\operatorname{spt}}

\newcommand{\clB}{B}       

\newcommand{\Lip}{\operatorname{Lip}}

\newcommand{\mass}{{\mathbf M}}

\newcommand{\AKc}{{\mathbf M}}    
\newcommand{\AKic}{{\mathbf I}}   

\begin{document}

\title[Plateau's problem in locally non-compact spaces]{Plateau's problem for integral currents in locally non-compact metric spaces}

\author{Stefan Wenger}

\address
  {Department of Mathematics\\
 University of Fribourg\\
  Chemin du Mus\'ee 23\\
  1700 Fribourg, Switzerland\\
  and\\
  Department of Mathematics\\
University of Illinois at Chicago\\
851 S. Morgan Street\\
Chicago, IL 60607--7045 }
\email{wenger@math.uic.edu {\rm or} stefan.wenger@unifr.ch}

\date{\today}

\thanks{Partially supported by NSF grants DMS 0956374 and DMS 1056263}


\begin{abstract}
 The purpose of this article is to prove existence of mass minimizing integral currents with prescribed possibly non-compact boundary in all dual Banach spaces and furthermore in certain spaces without linear structure, such as injective metric spaces and Hadamard spaces. We furthermore prove a weak$^*$-compactness theorem for integral currents in dual spaces of separable Banach spaces. Our theorems generalize results of Ambrosio-Kirchheim, Lang, the author, and recent results of Ambrosio-Schmidt.
\end{abstract}

\maketitle

\section{Introduction}\label{section:introduction}

In this paper we study the generalized Plateau problem in the context of locally non-compact metric spaces. Roughly speaking, this problem concerns the question of existence of an $(m+1)$-dimensional generalized surface of least volume with prescribed $m$-dimensional boundary in a given metric space. A suitable notion of surface in the context of area minimization problems is provided by the theory of integral currents. In the setting of Euclidean space, this theory was developed by Federer-Fleming in \cite{Federer-Fleming}, who solved Plateau's problem in the class of integral currents in Euclidean space. In \cite{Ambr-Kirch-curr}, Ambrosio-Kirchheim extended Federer-Fleming's theory to the setting of complete metric spaces. They then solved the generalized Plateau problem in the class of integral currents in compact metric spaces (provided that given a boundary the family of fillings is not empty) and furthermore in dual spaces of separable Banach spaces, provided that the prescribed boundary lies in a compact set. 
This result was generalized in \cite{Wenger-GAFA} to all dual Banach spaces as well as to Hadamard spaces (for Hadamard spaces, the result is due to U.~Lang), still requiring that the prescribed boundary remain in a compact set. It was only shown very recently by Ambrosio-Schmidt \cite{Ambrosio-Schmidt} that the condition on compact boundary can be dropped in the case of separable dual Banach spaces. 
%
%
The aim of the present paper is to remove the separability condition made in \cite{Ambrosio-Schmidt}. In fact, we will show that the generalized Plateau problem can be solved for non-compact boundaries in a class of metric spaces which includes e.g. all dual Banach spaces (also non-separable ones), all injective metric spaces, and all Hadamard spaces, i.e. complete simply-connected metric spaces of non-positive curvature in the sense of Alexandrov. 
We therefore generalize corresponding results in \cite{Ambr-Kirch-curr}, \cite{Ambrosio-Schmidt}, and \cite{Wenger-GAFA}. We furthermore partially generalize weak$^*$-compactness theorems proved in \cite{Ambr-Kirch-curr} and \cite{Ambrosio-Schmidt} to the setting of dual spaces of separable Banach spaces.

We now give precise formulations of our main results. Given a complete metric space $X$ and $m\geq 0$, we will denote by $\AKc_m(X)$ and $\AKic_m(X)$ the spaces of metric $m$-currents of finite mass and of integral $m$-currents, respectively, in the sense of Ambrosio-Kirchheim \cite{Ambr-Kirch-curr}. Given $T\in\AKc_m(X)$, the mass of $T$ will be denoted by $\mass(T)$ and, in case $m\geq 1$, the boundary of $T$ by $\bdry T$. We refer to \secref{section:prelims} for the basic definitions from the theory of metric currents. Our first theorem gives solutions to Plateau's problem and to the corresponding free boundary problem in the context of Banach spaces.

\bt\label{thm:banach}
 Let $X$ be a Banach space which is $1$-complemented in a dual Banach space. Let $m\geq 0$. Then 
 \begin{enumerate}
  \item for every $S\in\AKic_m(X)$ with $\bdry S=0$ there exists $T\in\AKic_{m+1}(X)$ with $\bdry T = S$ and such that
  \begin{equation*} \mass(T) \leq \mass(T')
 \end{equation*}
 for all  $T'\in\AKic_{m+1}(X)$ with $\bdry T'=S$;
  \item for every $S\in\AKc_m(X)$ there exists $T\in\AKic_{m+1}(X)$ such that $$\mass(T) +\mass(\bdry T - S) \leq \mass(T') +\mass(\bdry T' - S)$$ for all $T'\in\AKic_{m+1}(X)$.
 \end{enumerate}
\et

If $m=0$ then the condition $\bdry S=0$ in (i) should be replaced by $S(1)=0$, see \secref{section:prelims}. In the above, a Banach space $X$ is said to be $1$-complemented in a dual Banach space if $X$ is (isometric to) a subspace of a dual Banach space $Y$ such that there is a norm $1$ projection from $Y$ to $X$. Particular examples of such spaces are dual Banach spaces, $L^1$-spaces, and $L$-embedded Banach spaces.
Note that no compactness assumption on $\spt S$ or separability assumption on $X$ is made, and that furthermore $S$ in (ii) is only required to be a current of finite mass.
\thmref{thm:banach} generalizes \cite[Theorem 10.6]{Ambr-Kirch-curr}, \cite[Theorem 1.5]{Wenger-GAFA}, and Theorems 1.1 and 1.3 in \cite{Ambrosio-Schmidt}. 
Our next result provides a weak$^*$-compactness theorem which partially generalizes corresponding results in \cite{Ambr-Kirch-curr} and \cite{Ambrosio-Schmidt}.

\bt\label{thm:weakstar-cptness}
 Let $X$ be either a reflexive Banach space or the dual space of a separable Banach space. Let $m\geq 0$ and let $(T_n)\subset\AKic_m(X)$ be a sequence satisfying
  \begin{equation}\label{eq:unif-tot-mass-bd}
  \sup_n\left[\mass(T_n) + \mass(\bdry T_n)\right]<\infty
 \end{equation} 
and
 \begin{equation*}
  \lim_{r\to\infty}\left[\sup_n\|T_n\|(X\backslash B(0,r))\right]=0.
 \end{equation*} 
 Then there exists a subsequence $T_{n_j}$ which w$^*$-converges to some $T\in\AKic_m(X)$.
\et

Here, $B(0,r)$ denotes the (open) ball in $X$ of radius $r$ and center $0$. If $m=0$ then \eqref{eq:unif-tot-mass-bd} should be replaced by $\sup_n\mass(T_n)<\infty$. For the definition of w$^*$-convergence see \defref{def:weak-conv-var}. We remark that \thmref{thm:weakstar-cptness} fails in general if $X$ is the dual space of a non-separable Banach space, see Example~\ref{ex:non-separable-predual}. Our theorem generalizes \cite[Theorem 6.6]{Ambr-Kirch-curr} and \cite[Theorem 1.4]{Ambrosio-Schmidt} in the case of integral currents. While our theorem applies to a larger class of Banach spaces, \cite[Theorem 6.6]{Ambr-Kirch-curr} and \cite[Theorem 1.4]{Ambrosio-Schmidt} apply to normal currents as well. It is interesting to note however that \thmref{thm:weakstar-cptness} cannot be generalized to normal currents without the additional assumption that $X$ be separable, as was shown in \cite[Example B.1]{Ambrosio-Schmidt}.

We now briefly turn to metric spaces without a vector space structure. The methods which we use are not restricted to the setting of Banach spaces. Indeed, in \thmref{thm:1-complemented-metric-spaces} we will generalize \thmref{thm:banach} to a class of metric spaces. We will, in particular, obtain the following result as a consequence.
 
\bt\label{thm:Hadamard}
 Let $X$ be a Hadamard space or an injective metric space and let $m\geq 0$. Then assertions (i) and (ii) of \thmref{thm:banach} hold for $X$.
\et

Injective spaces are, by definition, absolute $1$-Lipschitz retracts. For the definition of Hadamard spaces we refer e.g. to \cite{Bridson-Haefliger}. 
\thmref{thm:Hadamard} generalizes a corresponding result for compact boundaries in Hadamard spaces which goes back to U.~Lang and which was published in \cite[Theorem 1.6]{Wenger-GAFA}. 

Finally, we mention that our approach differs from the one taken in \cite{Ambr-Kirch-curr} and \cite{Ambrosio-Schmidt}. Our methods combine arguments in the spirit of \cite{Wenger-GAFA} with a variant of a compactness theorem recently proved in \cite{Wenger-cpt} and \cite{Lang-Wenger-pointed}.

The paper is structured as follows. In \secref{section:prelims} we recall the definitions from the theory of currents needed for the paper. We furthermore define the notion of local weak convergence and establish a relationship with weak convergence. In \secref{section:proofs-mains} we prove a variant of the compactness theorems established in \cite{Wenger-cpt} and \cite{Lang-Wenger-pointed} and apply it in the proofs of the theorems stated above. We will furthermore state and prove some generalizations of the above results, see \thmref{thm:1-complemented-metric-spaces} and \thmref{thm:abstract-cptness}.

\section{Preliminaries}\label{section:prelims}

In this section we recall the basic definitions from the theory of metric currents developed in \cite{Ambr-Kirch-curr} which we will need in the sequel. With the exception of the definition of local weak convergence, its relationship with weak convergence, and \lemref{lem:conv-maps-currents}, all definitions and results appear in \cite{Ambr-Kirch-curr}.

Let $(X,d)$ be a complete metric space. We denote by $\Lip(X)$ and $\Lip_b(X)$ the spaces of real-valued Lipschitz functions and bounded Lipschitz functions on $X$, respectively. The Lipschitz constant of a Lipschitz function $f$ will be denoted by $\Lip(f)$. 
\bd\label{def:current}
Let $m\geq 0$. An $m$-dimensional metric current  $T$ on $X$ is a multi-linear functional $T:\Lip_b(X)\times\Lip^m(X)\to\R$ satisfying the following
properties:
\begin{enumerate}
 \item If $\pi^j_i\to \pi_i$ pointwise as $j\to\infty$ and if $\sup_{i,j}\Lip(\pi^j_i)<\infty$ then
       \begin{equation*}
         T(f,\pi^j_1,\dots,\pi^j_m) \longrightarrow T(f,\pi_1,\dots,\pi_m).
       \end{equation*}
 \item If $\{x\in X:f(x)\not=0\}$ is contained in the union $\bigcup_{i=1}^mB_i$ of Borel sets $B_i$ and if $\pi_i$ is constant 
       on $B_i$ then
       \begin{equation*}
         T(f,\pi_1,\dots,\pi_m)=0.
       \end{equation*}
 \item There exists a finite Borel measure $\mu$ on $X$ such that
       \begin{equation}\label{equation:mass-def}
        |T(f,\pi_1,\dots,\pi_m)|\leq \prod_{i=1}^m\Lip(\pi_i)\int_X|f|d\mu
       \end{equation}
       for all $(f,\pi_1,\dots,\pi_k)\in\Lip_b(X)\times\Lip^m(X)$.
\end{enumerate}
\ed
The space of $m$-dimensional metric currents on $X$ is denoted by $\AKc_m(X)$ and the minimal Borel measure $\mu$
satisfying \eqref{equation:mass-def} is called mass of $T$ and denoted by $\|T\|$. We also call mass of $T$ the number $\|T\|(X)$ 
which we denote by $\mass(T)$.
The support of $T$ is  the closed set $$\spt T = \{x\in X:  \text{ $\|T\|(B(x,r))>0$ for all $r>0$}\}.$$ 
As in \cite{Ambr-Kirch-curr} we will assume throughout this paper that the cardinality of any set is an Ulam number. This is consistent with the 
standard ZFC set theory. We then have that $\spt T$ is separable and furthermore that $\|T\|$ is concentrated on a 
$\sigma$-compact set, i.\ e.\ $\|T\|(X\backslash C) = 0$ for a $\sigma$-compact set $C\subset X$ (see \cite{Ambr-Kirch-curr}). This will be relevant in the proof of statement (ii) of \thmref{thm:banach}.

In the following we will often abbreviate $\pi=(\pi_1, \dots, \pi_m)$ and write $T(f,\pi)$ instead of $T(f,\pi_1,\dots, \pi_m)$. 
The restriction of $T\in\AKc_m(X)$ to a Borel set $A\subset X$ is given by 
\begin{equation*}
  (T\rstr A)(f,\pi):= T(f\chi_A,\pi).
\end{equation*}
This expression is well-defined since $T$ can be extended to a functional on tuples for which the first argument lies in 
$L^\infty(X,\|T\|)$. Furthermore, $T\rstr A\in\AKc_m(X)$ by Theorem 3.5 in \cite{Ambr-Kirch-curr}.
If $m\geq 1$ and $T\in\AKc_m(X)$ then the boundary of $T$ is the functional
\begin{equation*}
 \bdry T(f,\pi_1,\dots,\pi_{m-1}):= T(1,f,\pi_1,\dots,\pi_{m-1});
\end{equation*}
it satisfies conditions (i) and (ii) in \defref{def:current}. If it moreover satisfies (iii) in \defref{def:current} then $T$ is called a normal current.
By convention, elements of $\AKc_0(X)$ are also called normal currents. If $m=0$ and $T\in\AKc_0(X)$ then we say $\bdry T=0$ if $T(1)=0$. 
The push-forward of $T\in\AKc_m(X)$ 
under a Lipschitz map $\varphi$ from $X$ to another complete metric space $Y$ is given by
\begin{equation*}
 \varphi_\# T(g,\tau):= T(g\circ\varphi, \tau\circ\varphi)
\end{equation*}
for $(g,\tau)\in\Lip_b(Y)\times\Lip^m(Y)$. This defines a $m$-dimensional current on $Y$.
It follows directly from the definitions that $\bdry(\varphi_{\#}T) = \varphi_{\#}(\bdry T)$.

\bd\label{def:weak-conv-var}
 A sequence of currents $T_n\in\AKc_m(X)$ is said to converge weakly to $T\in\AKc_m(X)$ if 
 \begin{equation}\label{eq:ptwise-conv}
  T_n(f,\pi) \to T(f,\pi)
 \end{equation}
  for every $(f,\pi)\in \Lip_b(X)\times\Lip^m(X)$. If \eqref{eq:ptwise-conv} only holds for those $(f,\pi)$ for which $f$ also has bounded support, we say that $T_n$ converges locally weakly to $T$. If $X$ is a dual Banach space and \eqref{eq:ptwise-conv} holds for those $(f,\pi)$ for which the $f$ and $\pi_i$ are also weak$^*$-continuous, we say that $T_n$ w$^*$-converges to $T$.
  \ed
  
  It is clear that if $T_n$ converges locally weakly to $T$ and
 \begin{equation*}
  \lim_{r\to\infty}\left[\sup_n\|T_n\|(X\backslash B(x_0,r))\right]=0
 \end{equation*} 
 for some $x_0\in X$ then $T_n$ converges weakly to $T$. We furthermore have the following easy lemma.

 \bl\label{lem:conv-maps-currents}
  Let $X, Y$ be complete metric spaces, $m\geq 0$, and $S\in\AKc_m(X)$. Let $(\varphi_j)$ be a sequence of Lipschitz maps $\varphi_j:\spt S\to Y$  with $\sup_j\Lip(\varphi_j)<\infty$. If $\varphi_j$ converges pointwise to a Lipschitz map $\varphi:\spt S\to Y$ then $\varphi_{j\#}S$ converges weakly to $\varphi_{\#}S$.
 \el

\begin{proof}
 Let $(f,\pi)\in\Lip_b(Y)\times \Lip^m(Y)$. We may assume without loss of generality that $|f|\leq 1$ and that each $\pi_i$ is $1$-Lipschitz. Set $A:= \sup_j\Lip(\varphi_j)$ and let $\varepsilon>0$. By $\sigma$-compactness of $\spt S$ there exists $C\subset\spt S$ compact such that $\|S\|(X\backslash C)\leq \varepsilon/(4A^m)$. We then obtain for each $j$ that
 \begin{equation*}
  \begin{split}
   |\varphi_{\#}S(f,\pi) - \varphi_{j\#}S(f,\pi)| &\leq |(S\rstr C)(f\circ\varphi, \pi\circ\varphi) - (S\rstr C)(f\circ\varphi, \pi\circ\varphi_j)|\\
   &\quad\quad + |(S\rstr C)(f\circ\varphi - f\circ\varphi_j, \pi\circ\varphi_j)|+ 2A^m\|S\|(X\backslash C)\\
   &\leq |(S\rstr C)(f\circ\varphi, \pi\circ\varphi) - (S\rstr C)(f\circ\varphi, \pi\circ\varphi_j)|\\
   &\quad\quad +A^m\int_C|f\circ\varphi - f\circ\varphi_j|\,d\|S\| + \frac{\varepsilon}{2}.
  \end{split}
 \end{equation*}
 The first term after the last inequality sign converges to $0$ as $j\to\infty$ by the continuity property of currents, property (i) in \defref{def:current}. The second term converges to $0$ since $f\circ\varphi_j$ converges uniformly to $f\circ\varphi$ on the compact set $C$. Therefore, the sum of the three terms is $\leq \varepsilon$ for all $j$ large enough. This shows that $\varphi_{j\#}S$ converges weakly to $\varphi_{\#}S$.
\end{proof}

\sloppy
In this paper we will mainly be concerned with integral currents. 
\fussy
An element $T\in\AKc_0(X)$ is called integer rectifiable if there exist finitely many points $x_1,\dots,x_n\in X$ and $\theta_1,\dots,\theta_n\in\Z\backslash\{0\}$ such
that
\begin{equation*}
 T(f)=\sum_{i=1}^n\theta_if(x_i)
\end{equation*}
for all bounded Lipschitz functions $f$.
 A current $T\in\AKc_m(X)$ with $m\geq 1$ is said to be integer rectifiable if the following properties hold:
 \begin{enumerate}
  \item $\|T\|$ is concentrated on a countably $\hm^m$-rectifiable set and vanishes on all $\hm^m$-negligible Borel sets, where $\hm^m$ denotes the Hausdorff $m$-measure;
  \item For any Lipschitz map $\varphi:X\to\R^m$ and any open set $U\subset X$ there exists $\theta\in L^1(\R^m,\Z)$ such that 
    $\varphi_\#(T\rstr U)=\Lbrack\theta\Rbrack$, where $$\Lbrack\theta\Rbrack(f, \pi):= \int_{\R^m}\theta f\det\left(\nabla\pi\right)\,d\hm^k.$$
 \end{enumerate}
The space of integer rectifiable $m$-currents in $Z$ is denoted by $\AKic_m(Z)$.
Integer rectifiable normal currents are called integral currents. The corresponding space is denoted by $\AKic_m(X)$.

Recently, variants of Ambrosio-Kirchheim's theory that do not rely on the finite mass axiom have been developed by Lang in \cite{Lang-currents} and by Lang and the author in \cite{Lang-Wenger-pointed}. We will not need any definitions or results from these theories, however.
 
\section{Proofs of the main results}\label{section:proofs-mains}

The principal tool in the proofs of our main theorems is the following compactness result, which is a variant of the compactness theorem proven by the author in \cite{Wenger-cpt} and of a generalization proven by Lang and the author in \cite{Lang-Wenger-pointed}.

 \bl\label{lem:cptness-gen}
  Let $X$ be a complete metric space, $\Omega\subset X$ separable, and $m\geq 0$. Let $(T_n)\subset\AKic_m(X)$ be a sequence satisfying $$\sup_n[\mass(T_n) + \mass(\bdry T_n)]<\infty.$$ Then there exist a subsequence $(n_j)$, a complete metric space $Z$, and isometric embeddings $\varphi_j\colon X \hookrightarrow Z$ such that $\varphi_{j\#}T_{n_j}$ converges locally weakly to some $T\in\AKic_m(Z)$ and $\varphi_j|_\Omega$ converges pointwise to an isometric embedding $\varphi:\Omega\hookrightarrow Z$.
 \el

The statement of \lemref{lem:cptness-gen} can be strenghtened to apply also to the locally integral currents of \cite{Lang-Wenger-pointed} and local weak convergence can be replaced by convergence in the local flat topology. However, the above version suffices for all the applications in this paper. Note that, on the one hand, the results in \cite{Wenger-cpt} and \cite{Lang-Wenger-pointed} apply more generally to sequences $(T_n)$ of integral currents in a sequence $(X_n)$ of metric spaces. On the other hand, these results do not yield convergence of the isometric embeddings $\varphi_j$ as \lemref{lem:cptness-gen} does.

 The proof of \lemref{lem:cptness-gen} relies on the constructions used in \cite{Wenger-cpt} and \cite{Lang-Wenger-pointed}. Roughly speaking, the idea is to decompose each current $T_n$ into a sum $T_n = T_n^1 + T_n^2 + \dots$ of currents $T_n^i$ with mass growth $\geq \gamma_i r^m$ for some $\gamma_i>0$ independent of $n$. This mass growth implies that for fixed $i$ the sequence $(\spt T_n^i)$ of supports is a uniformly compact sequence of metric spaces and one can therefore use a variant of Gromov's compactness theorem together with the closure and compactness theorems for integral currents in a compact metric space \cite{Ambr-Kirch-curr} to produce a desired metric space $Z$ and a limit $T$ as in \lemref{lem:cptness-gen}. The decomposition procedure alluded above was proved in \cite[Lemma 5.1, Theorem 1.2]{Wenger-cpt}. It was summarized in Proposition 3.1 in \cite{Lang-Wenger-pointed} in a form which is suitable for our purposes.

In the proof below we will follow the proof of \cite[Theorem 1.1]{Lang-Wenger-pointed}. We would like to emphasize, however, that even though the paper \cite{Lang-Wenger-pointed} deals with locally integral currents we will only need integral currents in the sense of \cite{Ambr-Kirch-curr} and that we could instead follow the arguments given in the proof of \cite[Theorem 1.2]{Wenger-cpt}. 
 
 \begin{proof}
  Fix $x_0\in\Omega$ and choose numbers $0 < R_1< R_2 < \ldots \to\infty$ such that, after passing to a subsequence, we have $T_n\rstr\clB(x_0,R_r)\in\AKic_m(X)$ with
 \begin{equation*}
   \sup_n\mass(\bdry(T_n\rstr\clB(x_0,R_r))) <\infty
 \end{equation*}
for each $r\in\N$. Existence of such $R_r$ follows from \cite[Theorem 5.6]{Ambr-Kirch-curr} together with Fatou's Lemma. Set $R_0:= 0$, and define $A_{r}:= \clB(x_0,R_r)\setminus \clB(x_0,R_{r-1})$ and $$T_{r,n}:= T_n\rstr A_{r}$$
for $r\in\N$; clearly $T_{r,n}\in\AKic_m(X)$ and
\begin{equation*}
   \sup_n[\mass(T_{r,n})+ \mass(\bdry T_{r,n})] <\infty.
 \end{equation*}
Fix integers $1\leq j_1<j_2<\dots$ and positive numbers $\frac{1}{2}>\delta_1>\delta_2>\dots$ with $\sum_i\delta_i<\infty$. By possibly replacing $X$ by $\ell^\infty(X)$ we may assume that $X$ admits isoperimetric inequalities of Euclidean type for integral currents, see \cite[Corollary 1.3]{Wenger-GAFA}.
Now, let $$T_{r,n} = T_{r,n}^1+\dots+ T_{r,n}^{j_n+1}+U_{r,n}^1+ \dots+ U_{r,n}^{j_n+1}$$
be a decomposition with $T_{r,n}^1, \dots, T_{r,n}^{j_n+1},U_{r,n}^1, \dots, U_{r,n}^{j_n+1} \in\AKic_m(X)$ as in \cite[Propososition 3.1]{Lang-Wenger-pointed} for $T_{r,n}$, $R_r$, $x_0$, and $X$. Let $\{y_0, y_1, y_2, \dots \}\subset\Omega$ be a countable dense subset of $\Omega$ with $y_0=x_0$ and define $\Omega_s:= \{y_0,\dots, y_s\}$ for $s\in\N$. For $n,s\in\N$, define closed sets
\begin{equation*}
 B_n^s:= \Omega_s \cup  \bigcup_{r=1}^s\bigcup_{i=1}^{\min\{s, j_n\}}(\spt T_{r,n}^i \cup \spt U_{r,n}^i)
\end{equation*}
and note that $B_n^1 \subset B_n^2 \subset \ldots \subset X$. According to part~(i) of \cite[Propososition 3.1]{Lang-Wenger-pointed}, for each $s$, the sequence $(B_n^s)$ is uniformly compact. By~\cite[Proposition~5.2]{Wenger-cpt}, after passage to a subsequence, there exist isometric embeddings $\varphi_n\colon X\hookrightarrow Z$ and compact subsets $Y^1\subset Y^2\subset\ldots\subset Z$, for some complete metric space $Z$, such that 
\begin{equation*}
 \varphi_n(B_n^s)\subset Y^s
\end{equation*}
for all $n$ and $s$. It can be shown exactly as in the proof of \cite[Theorem 1.2]{Wenger-cpt} or the proof of \cite[Theorem 1.1]{Lang-Wenger-pointed} that, after passing to a further subsequence, $\varphi_{n\#}T_{n}$ converges locally weakly to some $T\in\AKic_m(Z)$. Note that the proof in \cite{Lang-Wenger-pointed} shows that after possibly replacing $Z$ by $\ell^\infty(Z)$ the $\varphi_{n\#}T_{n}$ converge to $T$ even in the local flat topology (which implies local weak convergence).

We are left to prove that, after passing to a further subsequence, $\varphi_n$ converges to an isometric embedding $\varphi:\Omega\hookrightarrow Z$.  For this, note first that since $\Omega_s\subset B_n^s$ we obtain $\varphi_n(\Omega_s)\subset Y^s$ for all $n$ and $s$. Since each $Y^s$ is compact we may assume, after passing to a subsequence, that $\varphi_n(y_s)$ converges to some $z_s\in Y^s$ as $n\to\infty$. By density of $\{y_0, y_1, \dots\}$ in $\Omega$ and by the fact that $d(z_s, z_r) = d(y_s, y_r)$ for all $r,s$ it follows that there exists an isometric embedding $\varphi: \Omega\hookrightarrow Z$ such that $\varphi(y_s) = z_s$ for all $s$. It finally follows that $\varphi_n|_\Omega$ converges to $\varphi$, which concludes the proof.
 \end{proof}
 
In the proofs of our main results we will use non-principal ultrafilters and ultralimits of sequences. Recall for this that a non-principal ultrafilter on $\N$ is a finitely additive probability measure $\omega$ on $\N$ (together with the $\sigma$-algebra of all subsets) such that $\omega$ takes values in $\{0,1\}$ only and
 $\omega(A)=0$ whenever $A\subset \N$ is finite.
Existence of non-principal ultra-filters on $\N$ follows from Zorn's lemma. It is not difficult to prove that if $(Y,\tau)$ is a compact Hausdorff topological space then for every sequence $(y_n)_{n\in\N}\subset Y$ there exists a unique point $y\in Y$ such that
\begin{equation*}
 \omega(\{n\in\N: y_n\in U\})=1
\end{equation*}
for every $U\in\tau$ containing $y$. We will call this point $y$ the ultralimit of the sequence $(y_n)$ and denote it by $\lim\nolimits_\omega y_n$.\\

We are now ready to prove Theorems~\ref{thm:banach} - \ref{thm:weakstar-cptness} stated in the introduction.

 \begin{proof}[Proof of \thmref{thm:banach}]
  We first prove (i). Let $S\in \AKic_m(X)$ with $\bdry S = 0$ and set $$s:= \inf\{\mass(T'): T'\in\AKic_{m+1}(X), \bdry T'=S\}.$$ Note that the set appearing on the right hand side is non-empty. Indeed, if $m=0$ then this follows from the fact that $X$ is geodesic. If $m\geq 1$ then this follows from the isoperimetric inequality \cite[Corollary 1.3]{Wenger-GAFA}. Let now $(T_n)\subset\AKic_{m+1}(X)$ be a sequence satisfying $\bdry T_n=S$ for all $n\in\N$ and such that $\mass(T_n)\to s$. Clearly, $$\sup_n[\mass(T_n) + \mass(\bdry T_n)]<\infty.$$ Set $\Omega:= \spt S$ and note that $\Omega$ is separable. Let $n_j$, $Z$, $\varphi_j$, $\varphi$, and $T$ be as in \lemref{lem:cptness-gen}, where $m$ is replaced by $m+1$. By \lemref{lem:conv-maps-currents}, $\varphi_{j\#}S$ converges weakly to $\varphi_{\#}S$ and hence $\bdry T = \varphi_{\#}S$. 
Now, view $X$ as a subspace of a dual Banach space $Y$ and let $\omega$ be a non-principal ultrafilter on $\N$. 
  We define a map $\psi: \spt T \to Y$ as follows. Let $z\in\spt T$. Since $\varphi_{j\#}T_{n_j}$ converges locally weakly to $T$ there exists a sequence $(x_j)\subset X$ with $x_j\in\spt T_{n_j}$ for all $j\in\N$ and such that $\varphi_j(x_j)\to z$. Clearly, $(x_j)$ is a bounded sequence in $Y$. 
  Since closed balls of finite radius in $Y$, endowed with the weak$^*$-topology, are compact and Hausdorff it follows that $(x_j)$ has an ultralimit $\lim_\omega x_j$ in $Y$.
  Define $\psi(z):= \lim_\omega x_j$. It follows from the lower semi-continuity of the norm in $Y$ with respect to weak$^*$-convergence that $\psi(z)$ is independent of the choice of sequence $(x_j)$ and that $\psi$ is $1$-Lipschitz. Since $\varphi_j(x) \to \varphi(x)$ for all $x\in\spt S$ it follows furthermore that $\psi\circ\varphi = \operatorname{id}_{\spt S}$. Finally, if $P: Y\to X$ is a projection of norm $1$ then $\hat{T}:= (P\circ \psi)_{\#}T$ satisfies $\hat{T}\in\AKic_{m+1}(X)$ and $$\bdry \hat{T} = (P\circ \psi)_{\#}(\varphi_{\#}S) = P_{\#} S = S$$ and $\mass(\hat{T})\leq \mass(T) \leq \liminf\mass(T_{n_j}) = s$. This completes the proof of (i).
 
 The proof of (ii) is analogous up to some minor modifications. Indeed, let $S\in\AKc_m(X)$ and set $$s:= \inf\{\mass(T')+\mass(\bdry T' - S): T'\in\AKic_{m+1}(X)\}.$$ Clearly, the set appearing on the right hand side is non-empty and thus $s$ is finite. Let $(T_n)\subset\AKic_{m+1}(X)$ be a sequence satisfying $\mass(T_n)+ \mass(\bdry T_n - S)\to s$. It follows that $$\sup_n[\mass(T_n) + \mass(\bdry T_n)]<\infty.$$ Set $\Omega:= \spt S$ and note that $\spt S$ is separable. Let $n_j$, $Z$, $\varphi_j$, $\varphi$, and $T$ be as in \lemref{lem:cptness-gen}, where $m$ is replaced by $m+1$. By \lemref{lem:conv-maps-currents}, $\varphi_{j\#}S$ converges weakly to $\varphi_{\#}S$. Let $Y$ and $P: Y\to X$ be as in (i). Define $\psi: \spt T\cup \spt S\to E$ in a similar way as $\psi$ was defined in (i). Then $\psi$ is $1$-Lipschitz and $\psi\circ\varphi = \operatorname{id}_{\spt S}$. It follows that $\hat{T}:= (P\circ \psi)_{\#}T$ satisfies $\hat{T}\in\AKic_{m+1}(X)$ and $$\bdry\hat{T} - S = (P\circ\psi)_{\#}(\bdry T - \varphi_{\#}S).$$ We thus obtain from the lower semi-continuity of mass with respect to (local) weak convergence that
 \begin{equation*}
  \mass(\hat{T}) + \mass(\bdry \hat{T} - S) \leq \mass(T) + \mass(\bdry T - \varphi_{\#}S)
  \leq \liminf_{j\to\infty}[\mass(T_{n_j}) + \mass(\bdry T_{n_j} - S)]
  =s.
  \end{equation*}
  This completes the proof of (ii).
 \end{proof}
 
 \br
We remark that the use of non-principal ultrafilters and ultralimits can be avoided in the above proof in the case that $X$ is the dual of a separable Banach space since in this case one can pass to a subsequence and use sequential weak$^*$-compactness of closed bounded balls, see also the proof of \thmref{thm:weakstar-cptness}.
\er

The proof of \thmref{thm:banach} easily generalizes to the following context. Let $(X, d)$ be a metric space and let $\omega$ be a non-principal ultrafilter on $\N$. A sequence $(x_n)_{n\in\N}\subset X$ is called bounded if $\sup_n d(x_n,x_1)<\infty.$ Define an equivalence relation $\sim$ on bounded sequences in $X$ by considering $(x_n)$ and $(x'_n)$ equivalent if $\lim_\omega d(x_n, x'_n) = 0$, and denote by $[(x_n)]$ the equivalence class of $(x_n)$ with respect to $\sim$. The ultra-completion $X_\omega$ of $X$ with respect to $\omega$ is the metric space given by the set $$X_\omega:= \left\{[(x_n)]: \text{ $(x_n)\subset X$ with $\sup d(x_n,x_1)<\infty$}\right\}$$ with the metric $d_\omega([(x_n)], [(x'_n)]):= \lim_\omega d(x_n, x'_n)$. We note that $X$ isometrically embeds into $X_\omega$ by the map which assigns to $x\in X$ the equivalence class of the constant sequence $(x)$. We may therefore view $X$ as a subset of $X_\omega$. In the following we will say that a metric space $X$ is $1$-complemented in some ultra-completion of $X$ if there exists a non-principal ultrafilter $\omega$ on $\N$ such that $X$, viewed as a subset of $X_\omega$, admits a $1$-Lipschitz retraction of $X_\omega$ onto $X$.
Note that if $X$ is a dual Banach space then $X$ is $1$-complemented in every ultra-completion of $X$. Consequently, if $X$ is a Banach space which is $1$-complemented in some dual Banach space (in the terminology established before the statement of \thmref{thm:banach}) then $X$ is $1$-complemented in every ultra-completion of $X$. The following result generalizes \thmref{thm:banach}.

\bt\label{thm:1-complemented-metric-spaces}
 Let $X$ be a complete metric space and suppose that $X$ is $1$-complemented in some ultra-completion of $X$. Let $m\geq 0$. Then 
 \begin{enumerate}
  \item for every $V\in\AKic_{m+1}(X)$ there exists $T\in\AKic_{m+1}(X)$ with $\bdry T = \bdry V$ and
  \begin{equation*} 
  \mass(T) \leq \mass(T')
 \end{equation*}
 for all  $T'\in\AKic_{m+1}(X)$ with $\bdry T'=\bdry V$;
  \item for every $S\in\AKc_m(X)$ there exists $T\in\AKic_{m+1}(X)$ such that $$\mass(T) +\mass(\bdry T - S) \leq \mass(T') +\mass(\bdry T' - S)$$ for all $T'\in\AKic_{m+1}(X)$.
  \end{enumerate}
\et

Note that the statement in (i) of \thmref{thm:1-complemented-metric-spaces} is slightly weaker than the statement (i) in \thmref{thm:banach} inasmuch as we assume that there exists a filling of $\bdry V$. This is needed since in the generality considered in \thmref{thm:1-complemented-metric-spaces} integral currents without boundary need not have a filling.

\begin{proof}
 The proof is analogous to the proof of \thmref{thm:banach}.
\end{proof}

\thmref{thm:Hadamard} now comes as a consequence of \thmref{thm:1-complemented-metric-spaces}. Indeed, if $X$ is a Hadamard space then so is every ultra-completion $X_\omega$ of $X$ and, as a closed convex subspace of $X_\omega$, the nearest point projection from $X_\omega$ to $X$ is $1$-Lipschitz, see e.g. \cite[Proposition II.2.4]{Bridson-Haefliger}. Therefore, in the terminology established above, $X$ is $1$-complemented in every ultra-completion of $X$. Statement (ii) of \thmref{thm:Hadamard} therefore follows from statement (ii) of \thmref{thm:1-complemented-metric-spaces}. As for (i), if $S\in\AKic_m(X)$ is such that $\bdry S = 0$ then by \cite[Corollary 1.4]{Wenger-GAFA} there exists $V\in\AKic_{m+1}(X)$ with $\bdry V = S$. Statement (i) of \thmref{thm:Hadamard} now follows from this together with statement (i) of \thmref{thm:1-complemented-metric-spaces}. If $X$ is an injective space then $X$ is an absolute $1$-Lipschitz retract and thus there exists a $1$-Lipschitz projection of $X_\omega$ onto $X$ and thus $X$ is $1$-complemented in every ultra-completion of $X$. Statements (i) and (ii) now follow as above.

We finally turn to the proof of our w$^*$-compactness theorem.

\begin{proof}[Proof of \thmref{thm:weakstar-cptness}]
 Let $(T_n)$ be as in the statement of the theorem. Set $\Omega:= \{0\}$ and let $n_j$, $Z$, $\varphi_j$, $\varphi$, and $T$ be as in \lemref{lem:cptness-gen}. Then $\varphi_j(0)$ converges to some $z^0\in Z$. Let $(z^k)\subset\spt T$ be a dense sequence. For each $k\in\N$ choose a sequence $(x_j^k)\subset X$ with $x_j^k\in\spt T_{n_j}$ for all $j$ and such that $\varphi_j(x_j^k)\to z^k$. Since $$\|x_j^k\| = d(\varphi_j(x_j^k), \varphi_j(0))\to d(z^k, z^0)$$ as $j\to\infty$ it follows that for $k$ fixed the sequence $(x_j^k)$ is bounded. After passing to a subsequence we may therefore assume that $x_j^k$ converges to some $x^k$ in the weak$^*$-topology of $X$. (Note that for this separability of a predual or reflexivity of $X$ is needed.) Define $\psi(z^k):= x^k$ and note that $\psi$ is $1$-Lipschitz by the lower semi-continuity of the norm on $X$ with respect to w$^*$-convergent sequences. Since the sequence $(z^k)$ is dense in $\spt T$ we can extend $\psi$ to a $1$-Lipschitz map $\psi:\spt T\to X$. Set $\hat{T}:=\psi_{\#}T$ and note that $\hat{T}\in\AKic_m(X)$. We will show that $T_{n_j}$ is w$^*$-convergent to $\hat{T}$. For this let $f, \pi_1, \dots,\pi_m\in\Lip(X)$ be weak$^*$-continuous and such that $|f|\leq C$ for some $C$. For each $j$ let $f^j:Z\to \R$ be a Lipschitz extension of $f\circ\varphi_j^{-1}$ with $\Lip(f^j)=\Lip(f)$ and such that $|f^j|\leq C$. Similarly, let $\pi_i^j:Z\to \R$ be a Lipschitz extension of $\pi_i\circ\varphi_j^{-1}$ with $\Lip(\pi_i^j)=\Lip(\pi_i)$. For each $k$ we have 
 \begin{equation*}
 \begin{split}
  |f\circ\psi(z^k) - f^j(z^k)| &\leq |f(x^k) - f(x_j^k)| + |f^j(\varphi_j(x_j^k)) - f^j(z^k)|\\
  & \leq |f(x^k) - f(x_j^k)| + \Lip(f)\cdot d(\varphi_j(x_j^k), z^k),
  \end{split}
 \end{equation*}
 from which it follows together with the weak$^*$-continuity of $f$ that $f^j$ converges pointwise to $f\circ\psi$ on $\spt T$. Analogously, $\pi_i^j$ converges to $\pi_i\circ\psi$ on $\spt T$. 
 Now, since $\varphi_{j\#}T_{n_j}$ converges locally weakly to $T$ and
  \begin{equation*}
  \lim_{r\to\infty}\left[\sup_n\|T_n\|(X\backslash B(0,r))\right]=0
 \end{equation*} 
it follows that $\varphi_{j\#}T_{n_j}$ converges weakly to $T$. By possibly replacing $Z$ by $\ell^\infty(Z)$ we may assume by \cite[Theorem 1.4]{Wenger-flatconv} that $\varphi_{j\#}T_{n_j}$ converges even in the flat norm to $T$. In particular, there exist $U_j\in\AKic_m(Z)$ and $V_j\in\AKic_{m+1}(Z)$ such that $T - \varphi_{j\#}T_{n_j} = U_j + \bdry V_j$ and $\mass(U_j) + \mass(V_j)\to0$. We now obtain
 \begin{equation*}
  \begin{split}
   |\hat{T}(f,\pi) - T_{n_j}(f,\pi)| &\leq |T(f\circ\psi, \pi\circ\psi) - T(f^j,\pi^j)| + |T(f^j, \pi^j) - \varphi_{j\#}T_{n_j}(f^j, \pi^j)|\\
    &\leq |T(f\circ\psi, \pi\circ\psi) - T(f\circ\psi,\pi^j)| + |T(f\circ\psi - f^j, \pi^j)|\\
    &\quad + |U_j(f^j, \pi^j)| + |V_j(1,f^j, \pi^j)|\\
    &\leq |T(f\circ\psi, \pi\circ\psi) - T(f\circ\psi,\pi^j)| + \prod_{i=1}^m\Lip(\pi_i) \int_Z|f\circ\psi - f^j|\,d\|T\|\\
    &\quad + \prod_{i=1}^m\Lip(\pi_i)\left[C + \Lip(f)\right]\cdot(\mass(U_j)+\mass(V_j)).
  \end{split}
 \end{equation*}
 Since each $\pi_i^j$ converges pointwise to $\pi_i\circ\psi$ on $\spt T$ with bounded Lipschitz constants, the first term after the last inequality sign converges to $0$ by the continuity property of currents. Since $f^j$ converges uniformly to $f\circ\psi$ on compact subsets of $\spt T$ and $\spt T$ is $\sigma$-compact, it follows that the second term converges to $0$ as well. Since also the third term converges to $0$ it follows that $T_{n_j}$  indeed w$^*$-converges to $\hat{T}$ as claimed. This concludes the proof.
 \end{proof}

The following example shows that in general the assumption in the statement of \thmref{thm:weakstar-cptness} that a predual of $X$ be separable cannot be dropped.

\begin{example}\label{ex:non-separable-predual}
 Let $X$ be the dual space of $\ell^\infty$. For $n\geq 1$ define $x_n\in X$ by $x_n(a):= a_n$ for $a=(a_1,a_2, \dots)\in\ell^\infty$, and let
 $T_n\in\AKic_0(X)$ be given by $T_n(f):= f(x_n)$ for every $f\in\Lip_b(X)$. It follows that $T_n$ is supported in the closed unit ball in $X$ and that $\mass(T_n) =1$ for every $n$. It is easy to show that there cannot exist a subsequence $(n_j)$ such that $T_{n_j}$ is w$^*$-convergent to some $T$. Indeed, given a subsequence $(n_j)$, let $a=(a_1,a_2, \dots)\in\ell^\infty$ be defined by $a_k = 1$ if $k=n_j$ for some $j$ even and $a_k=0$ otherwise. Define a weak$^*$-continuous function $f\in\Lip_b(X)$ by $f(x):= \varphi(x(a))$, where $\varphi$ is the truncation function given by $\varphi(t)= \max\{-1, \min\{1, t\}\}$ for $t\in\R$. Clearly, $T_{n_j}(f) = f(x_{n_j})$ does not converge as $j\to\infty$.
\end{example}

We note that the only facts about dual spaces $X$ of separable Banach spaces which are used in the proof of \thmref{thm:weakstar-cptness} are the sequential weak$^*$-compactness of closed bounded balls in $X$ and the lower semi-continuity of the norm on $X$ with respect to weak$^*$-convergent sequences. The proof therefore easily gives the following generalization of \thmref{thm:weakstar-cptness}.

\bt\label{thm:abstract-cptness}
 Let $X$ be a Banach space and $\mathcal{W}$ a vector space topology on $X$ such that closed bounded balls in $X$ are sequentially $\mathcal{W}$-compact and such that the norm on $X$ is lower semi-continuous with respect to $\mathcal{W}$-convergent sequences. Then the conclusion of \thmref{thm:weakstar-cptness} holds when w$^*$-convergence is replaced by $\mathcal{W}$-convergence.
\et

We note that $\mathcal{W}$-convergence of $T_{n_j}$ to $T$ means by definition that $T_{n_j}(f,\pi)\to T(f,\pi)$ for all $f, \pi_1,\dots, \pi_m\in\Lip(X)$ which are also $\mathcal{W}$-continuous and such that $f$ is bounded. \thmref{thm:abstract-cptness} generalizes Theorem A.1 of \cite{Ambrosio-Schmidt} for integral currents.

\end{document}